\documentclass{article}
\usepackage{amssymb}
\usepackage{amsmath}

\usepackage{cite}
\usepackage{amssymb}
\usepackage{amsfonts}
\usepackage{amscd}
\usepackage[matrix]{xypic}
\newtheorem{theorem}{Theorem}

\newtheorem{definition}[theorem]{Definition}

\newtheorem{lemma}[theorem]{Lemma}

\newtheorem{proposition}[theorem]{Proposition}

\newcommand{\R}{\Bbb{R}}
\newcommand{\C}{\Bbb{C}}
\newcommand{\g}{\frak{g}}
\newcommand{\A}{\mathcal{A}}
\newcommand{\K}{\Bbb{K}}

\begin{document}

\title{\bf Rigid current Lie algebras}

\author{Michel
GOZE \thanks{M.Goze@uha.fr.} - 
Elisabeth REMM \thanks{%
corresponding author: e-mail: E.Remm@uha.fr} \\
\\
{\small Universit\'{e} de Haute Alsace, F.S.T.}\\
{\small 4, rue des Fr\`{e}res Lumi\`{e}re - 68093 MULHOUSE - France}}
\date{}
\maketitle

\noindent {\bf 2000 Mathematics Subject Classification.} Primary  17Bxx, Secondary 16Bxx.

\bigskip

\noindent {\bf Keywords.} 
Current Lie algebras, tensor product of algebras, rigidity, deformations of current Lie algebras.

\begin{abstract}
A current Lie algebra is contructed from a tensor product of a Lie algebra and a commutative associative 
algebra of dimension greater than $2$.
In this work we are interested in deformations of such algebras and in the problem of rigidity.
In particular we prove that a current Lie algebra is rigid if it is isomorphic 
to a direct product $\g\times \g\times ...\times \g$ where $\g$ is a rigid Lie algebra.
\end{abstract}

\section{Current Lie algebras}

If $\g$ is a Lie algebra over a  field $\K$ and $ \mathcal{A}$ a $\K$- associative commutative algebra, then
$\g \otimes \mathcal{A}$, provided with the bracket
$$[X\otimes a,Y\otimes b]=[X,Y]\otimes ab$$
for every $X,Y \in \g$ and $a,b \in \mathcal{A}$ is a Lie algebra. If $dim (\mathcal{A})=1$ such an algebra is isomorphic to
$\g$. If $dim (\mathcal{A})>1$ we will say that $\g \otimes \mathcal{A}$ with the previous bracket is a current Lie algebra.

In \cite{R.G.prodtens} we have shown that if  $\mathcal{P}$ is a quadratic operad , there is an associated quadratic operad, noted 
 $\tilde{\mathcal{P}}$ such that the tensor product of a  $\mathcal{P}$-algebra by a  $\tilde{\mathcal{P}}$-algebra
is a  $\mathcal{P}$-algebra for the natural product. In particular, 
if the operad $\mathcal{P}$ is $\mathcal{L}ie$, then $\tilde {\mathcal{L}ie} =\mathcal{L}ie^!=
\mathcal{C}om$ and a $\mathcal{C}om$-algebra is a commutative associative algebra. In this context we find again the 
notion of current Lie algebra.

In this work we study the deformations of a current Lie algebra and we show that a current Lie algebra is rigid if and only if it is isomorphic to
$\g \times \g \times ...\times \g$ where $\g$ is a rigid Lie algebra. The notion of rigidity is related to the second group of the Chevalley cohomology.
For the current Lie algebras, this group is not wellknown. Recently some relation between 
$H^2(\g \otimes \mathcal{A},\g \otimes \mathcal{A})$ and $H^2(\g,\g)$ and $H^2_H(\mathcal{A},\mathcal{A})$
are given in \cite{Z} but often when $\g$ is abelian. Let us note also that the scalar cohomology 
has been studied in \cite{NW}.

\section{Determination of rigid current  Lie algebras}

\subsection{On the rigidity of Lie algebras}

Let us remind briefly some properties of the variety of Lie  algebras (for more details, see \cite{A}). Let
$\g$ be a $n$-dimensional $\K$-Lie algebra. Since the underlying vector space is isomorphic to $\K^n$, there exists a one to one
correspondance between the set of Lie brackets of $n$-dimensional Lie algebras and the skew-symmetric bilinear maps
$\mu : \K^n \times \K^n\rightarrow \K^n$ satisfying the Jacobi identity. We  denote by $\mu _{\g}$ this bilinear map 
corresponding to $\g$. In this framework, we can identify $\g$ with the pair $(\K^n,\mu_{\g})$. Let us fix definitively a basis
$\{X_1,...,X_n\}$ of $\K^n$. The structure constants $(C_{ij}^k)$ of $\mu _{\g}$ are given by
$$\mu_{\g}(X_i,X_j)=\sum _{k=1}^nC_{ij}^k \ X_k$$
and we can identify $\mu_{\g}$ with the $N$-uple $(C_{ij}^k)$ with $N=\frac{n^2(n-1)}{2}$. The Jacobi identity
satisfied by $\mu_{\g}$ is equivalent to the polynomial system :
\begin{eqnarray}
\label{Jacobi}
\sum_{l=1,...,n} C_{ij}^lC_{lk}^s+C_{jk}^lC_{li}^s+C_{ki}^lC_{lj}^s=0.
\end{eqnarray}
Thus a Lie algebra is a point of $\K^N$ whose  coordinates $(C_{ij}^k)$ satisfy (\ref{Jacobi}). 
So the set of $n$-dimensional Lie algebras on $\K$ is identified with the algebraic variety $L_n$ embedded into
$\K^N$ and defined by the system of polynomial equations (\ref{Jacobi}). We will always denote by $\mu$ a point
of $L_n$. The algebraic group $GL(n,\K)$ acts on $L_n$ by:
\begin{eqnarray}
\label{action}
(f,\mu) \in GL(n,\K)\times L_n\longrightarrow \mu_f \ \in L_n
\end{eqnarray}
where $\mu_f$ is given by $\mu_f(X,Y)=f^{-1}(\mu(f(X),f(Y))$ for every $X,Y \in \K ^n.$ The orbit $\mathcal{O}(\mu)$
of $\mu$ related to this action corresponds to the  Lie algebras isomorphic to $\g=(\mu,\K^n).$ We provide the algebraic variety
$L^n$ with the Zariski topology.
\begin{definition}{\label{rigide}}
The Lie algebra $\g=(\mu,\K^n)$ is rigid if the orbit $\mathcal{O}(\mu)$ is open in $L_n$.
\end{definition}
Let us suppose that $\K$ is an algebraically closed field. A way of building rigid Lie algebras rests on the 
Nijenhuis Richardson Theorem : Let $H^*(\g,\g)$ be the Chevalley cohomology of $\g$. If $H^2(\g,\g)=0$ then
$\g$ is rigid. Let us note that the converse is false, numerous examples are described in \cite{A} (in fact a rigid Lie algebra
whose cohomology $H^2(\g,\g)$ is not trivial is such that the affine schema $\mathcal{L}_n$ given by the Jacobi
ideal is not reduced to the point $\mu$ defining $\g$.)

\noindent An intuitive way of defining the notion of rigidity is to consider a rigid algebra as indeformable, 
that is any close algebra is isomorphic to it. A general definition of deformations was proposed in \cite{G.R.valued}.
Let $A$ be a commutative $\K$ algebra of valuation such that the residual field $A/\frak{m}$ is isomorphic to $\K$
where $\frak{m}$ is the maximal ideal of $A$. If $\g$ is a $\K$-Lie algebra then the tensor product 
$\g \otimes A$ is an $A$-algebra denoted by $\g _A$.
\begin{definition}{\label{deformation}}
A deformation of $\g$ is an $A$-Lie algebra $\g'_A$ such that the underlying $A$-module is $\g_A$ and the brackets
$[u,v]_{\g'_A}$ and $[u,v]_{\g_A}$ of $\g'_A$ and $\g_A$ satisfy 
$$\lbrack u,v]_{\g'_A}-[u,v]_{\g_A} \in \g \otimes \frak{m}.$$
\end{definition}
When $A=\C[[t]]$ we find the classical notion of deformation given by Gerstenhaber. When $A$ is the ring of limited elements
in a Robinson nonarchimedean extension of $\C$, we find the notion of perturbations \cite{G}. If $\g'_A$ is a deformation of $\g$
then we have
$$\lbrack u,v]_{\g'_A}-[u,v]_{\g_A} =\sum_{i=1}^k \ \epsilon _1\epsilon _2...\epsilon _i\phi _i$$
where $\epsilon _i \in \frak{m}$ and $\{\phi _1,...,\phi _k\}$ a family of 
independant skewsymmetric bilinear maps on $\K^n\times \K^n$ with values in $\K^n$. In particular $\phi_1 \in Z^2(\g,\g)$
and if $\g'_A$ is isomorphic to $\g_A$ this map belongs to $B^2(\g,\g)$. We deduce that the deformations of $\g$
are parametrized by $H^2(\g,\g).$
\noindent In the following, we are going to determine the current Lie algebras which are rigid.

\subsection{The manifold $L_{(p,q)}$}

Let $\frak{g}=\frak{g} =\frak{g}_p \otimes \mathcal{A}_q$ be a current Lie algebra  where $\frak{g}_p$
 is a $p$-dimensional Lie algebra and $\mathcal{A}_q$ a $q$-dimensional associative commutative algebra.
We suppose that  $\mathbb{K}$ is an algebraically closed field  of  characteristic $0$.
Let  $\{ X_1 ,..., X_p\}$ be a basis of $\frak{g}_p$ and $\{ e_1 ,..., e_p\}$ a basis of $\mathcal{A}_q$. 
If we denote by $\{C_{ij}^k \}$ and $\{D_{ab}^c \}$ the structure constants of  $\frak{g}_p$ and $\mathcal{A}_q$ 
with regards to these basis, then the Lie bracket $\mu_{\g}=\mu_{\g_p}\otimes \mu_{\mathcal{A}_q}$ of $\g$ where
$\mu_{\g_p}$ is the multiplication of $\g$ and $\mu_{\mathcal{A}_q}$ the multiplication of 
$\mathcal{A}_q$, satisfy:
$$\mu_{\g}(X_i \otimes e_a, X_j \otimes e_b)=\sum_{k,c} C_{ij}^k D_{ab}^c X_k \otimes e_c,$$
and the structure constants of $\frak{g}$ with respect to the basis $\{ X_i \otimes e_a \}_{i=1,...,p; \ a=1,...,q}$ are $\{ C_{ij}^k D_{ab}^c \}.$
The Jacobi relations so are written
$$\sum_{l,r} C_{ij}^l C_{lk}^s D_{ab}^r D_{rc}^t +C_{jk}^l C_{li}^s D_{bc}^r D_{ra}^t +C_{ki}^l C_{lj}^s D_{ca}^r D_{rb}^t=0$$
$\forall (s,t) \in \left\{ \{1,...,p \} \times \{1,...,q \} \right\}.$
These polynomial relations define a structure of algebraic variety denoted by $L_{(p,q)}$ and embedded
in the vector space whose coordinates are the structure constants  $\{C_{ij}^k D_{ab}^c  \}$. It is a closed subvariety
of  $L_{pq}$. Let $G(p,q)$ be the algebraic group $G(p,q)=GL(p) \otimes GL(q).$ This group acts naturally
on $L_{(p,q)}$ by
$$(f \otimes g).(\mu_{\g_p}\otimes \mu_{\mathcal{A}_q})(X\otimes a,Y \otimes b)=
f^{-1}(\mu_{\g_p}(f(X),f(Y)))\otimes g^{-1}(\mu_{\mathcal{A}_q}(g( a),g( b))).$$ 
We denote by $\mathcal{O}_{p,q}(\frak{g}_p \otimes \mathcal{A}_q )$ the orbit in $L_{(p,q)}$ of $\mu_{g}$ 
corresponding to this action.

\smallskip

\noindent Thus there are two types of deformations:

- The deformations of $\frak{g}$ in the manifold $L_{pq}$.
These deformations are parametrized by the second Chevalley cohomology space $H^2(\frak{g},\frak{g})$.

- The deformations of $\frak{g}$ in the manifold $L_{(p,q)}$. They are parametrized by 
$H^2_C(\frak{g}_p,\frak{g}_p) \oplus H^2_H(\mathcal{A}_q,\mathcal{A}_q)$ where $H^2_H(\mathcal{A}_q,\mathcal{A}_q)$  
is the Harrison cohomology of the associative commutative algebra $\mathcal{A}_q$. 
\begin{definition}
The Lie algebra $\frak{g}_p \otimes \mathcal{A}_q$ is rigid in $L_{(p,q)}$ if the orbit 
$\mathcal{O}_{p,q}(\mu_{\g})$ is open (in the Zariski sense). It is rigid if the orbit $\mathcal{O}(\mu_{\g})$
 related to the 
action of $GL(pq)$ in $L_{pq}$ is open. 
\end{definition}
It is clear that the rigidity implies the rigidity in $L_{(p,q)}$.

\begin{proposition}
A current Lie algebra  $\g=\frak{g}_p \otimes \mathcal{A}_q$ is rigid in $L_{(p,q)}$
if and only if $\frak{g}_p$ is rigid in $L_{p}$ and $\mathcal{A}_q$ is rigid in $\mathcal{C}om(q)$,
the variety of $q$-dimensional associative commutative $\K$ -algebras.
\end{proposition}

\noindent {\bf Example.} $p=2, q=2$ $(\mathbb{K}=\mathbb{C})$
There is, up to isomorphisms, only one $2$-dimensional rigid Lie algebra . It is defined by 
$[X_1,X_2]=X_2$. There is only one $2$-dimensional associative commutative algebra. 
 It is given by $e_1^2=e_1, e_2^2=e_2,e_1e_2=0$ and corresponds to the semi-simple algebra
$A_1^2=M_1(\mathbb{K})\times M_1(\mathbb{K})$ where $M_n(\mathbb{K})$ is the algebra of $n$-matrices on $\mathbb{K}.$
The Lie algebra $\frak{g}_2\otimes A_1^2$ is rigid in $L_{(2,2)}$. This algebra is isomorphic to 
$\frak{g}_2 \times \frak{g}_2.$ It is also rigid in $L_4.$

\subsection{Structure of  rigid  current Lie algebras}

Recall that if $\frak{g}$ is a finite dimensional rigid Lie algebra it admits a decompostion 
$\frak{g}=\frak{s} \oplus \frak{t} \oplus \frak{n}$ where
$\frak{t} \oplus \frak{n}$ is the radical of $\frak{g}$, $\frak{t}$ is a maximal abelian subalgebra 
whose adjoint operators 
$ad\, X , X \in \frak{t}$ are semi-simple and $\frak{n}$ is the nilradical.
If $\frak{g}=\frak{g}_p \otimes \mathcal{A}_q$ is rigid then  $\frak{g}_p$ is rigid in $L_p$.
If $\frak{g}_p$ is solvable then $\frak{g}$ too and we have 
$$\frak{g}_p=\frak{t}_p \oplus \frak{n}_p \ \mbox{ \rm{ and}}  \ \ \frak{g}=\frak{t} \oplus \frak{n}.$$
Since $\frak{n}_p \otimes \mathcal{A}_q$ is a nilpotent ideal of $\frak{g},$ 
$\frak{n}_p \otimes \mathcal{A}_p \subset \frak{n}.$

\begin{lemma}
If $\frak{g}=\frak{g}_p \otimes \mathcal{A}_q$ is rigid, then $\mathcal{A}_q$ has a non zero idempotent.
\end{lemma}

\noindent {\it Proof.} If $\mathcal{A}_q$ is a nilalgebra then $\frak{g}$ is nilpotent. In fact
if $X \in \frak{g}_p$ and $a \in \mathcal{A}_q$ we have $[ad(X \otimes a)]^m=(ad\, X)^m \otimes (L_a)^m$
where $L_a:\mathcal{A}_q \rightarrow \mathcal{A}_q$ is the left multiplication by $a$. Since $\mathcal{A}_q$
is a nilalgebra, every element is nilpotent and there exits $m_0$ such that $(L_a)^{m_0}=0.$
Thus $ad ( X\otimes a)$ is a nilpotent operator for any $X$ and $a$. 
This implies that $\frak{g}$ is nilpotent.
Let $f$ be a derivation of $\frak{g}_p.$ Then $f\otimes Id$ is a derivation of $\frak{g}.$ Since 
$\frak{g}_p$ is rigid, we can find a inner non trivial derivation $ad \, X$ which is diagonal. In this case 
$ad \, X \otimes Id$ is a non trivial diagonal derivation of $\frak{g}$. By hypothesis
 $\frak{g}$ is rigid. But any rigid nilpotent Lie algebra is characteristically nilpotent, that is, 
every derivation is nilpotent. We have a contradiction and  $\mathcal{A}_p$ can not be
a nilalgebra. Since it is finite dimensional, it admits a non zero idempotent.

\begin{proposition}
If $\frak{g}=\frak{g}_p \otimes \mathcal{A}_q$ is rigid then $\mathcal{A}_q$ is an associative commutative
rigid unitary algebra in $\mathcal{C}om(q)$.
\end{proposition}

\noindent {\it Proof.} Let $e$ be in $\mathcal{A}_q$ and satisfying $e^2=e.$ The associated Pierce decomposition
$$\mathcal{A}_q= \mathcal{A}^{00}_q \oplus \mathcal{A}^{10}_q \oplus \mathcal{A}^{01}_q \oplus \mathcal{A}^{11}_q$$ 
where 
$$\mathcal{A}^{ij}_q=\{ x \in \mathcal{A}_q \ \mbox{\rm \ such that} \ \  e\cdot x=i x,x \cdot e=jx \}$$
 reduce to 
$\mathcal{A}_q=\mathcal{A}^{11}_q \oplus \mathcal{A}^{00}_q  $ because $\mathcal{A}_q $ is commutative and we have
$\mathcal{A}^{11}_q \cdot \mathcal{A}^{00}_q =\{ 0 \}.$
Thus $\mathcal{A}_q$ is a direct sum of two commutative algebras. Since $\mathcal{A}_q $  is rigid, the algebras
 $\mathcal{A}^{11}_q$ and $\mathcal{A}^{00}_q $ are also rigid. The subalgebra $\mathcal{A}^{11}_q $ is unitary 
($e$ is the unit element). From the previous lemma $\mathcal{A}^{00}_q $ has an idempotent and admits a decomposition
$$\mathcal{A}^{00}_q=\mathcal{A}^{0011}_q \oplus \mathcal{A}^{0000}_q $$ 
with $\mathcal{A}^{0011}_q \neq \{ 0 \}.$ By induction we deduce that 
$$\mathcal{A}_q =\mathcal{A}^{1}_q \oplus ... \oplus \mathcal{A}^{p}_q $$
with $\mathcal{A}^{i}_q $ with unit $e_i$ and $\{ e_1,...,e_p \}$ is a system of pairwise orthogonal idempotents. 
Then $e_1 + ...+e_p$ is a unit of $\mathcal{A}_q .$
 
\begin{theorem}
Let $\frak{g}_p$ be a rigid Lie algebra with solvable non nilpotent radical such that $Z(\frak{g})=\{ 0 \}.$
Then $\frak{g}=\frak{g}_p \otimes \mathcal{A}_q$ is rigid if and only if $\mathcal{A}_q=M_1^q$ is given by
$$e_i^2=e_i \ , i=1,...,q \ \  \mbox{\rm and} \ \ e_i \cdot e_j=0 \ \mbox{\rm if} \ i \neq j.$$ 
\end{theorem}

\noindent {\it Proof.} Since $\mathcal{A}_q$ is unitary, the radical of $\frak{g}$ solvable and non nilpotent. 
Moreover $Z(\frak{g}_p)=\{ 0 \}$ implies that $Z(\frak{g})=\{ 0 \}.$ In fact if 
$U= \sum_{j,a} \alpha_{ja} X_j \otimes x_a$ is in the center of $\frak{g}$, then $[U,X \otimes 1]=0$ for each
$X \in \frak{g}_p.$ Thus 
$$\sum \alpha_{j,a}[X_j ,X] \otimes x_a=0.$$
We have  $[\sum_j \alpha_{ja}X_j,X]=0$ for each $a$ and $X.$ 
So $\sum_j \alpha_{ja}X_j \in Z(\frak{g}_p)$ for any $a.$ Therefore $\alpha_{ja}=0$ for any $a$ and $U=0.$

Consequently $\frak{g}$ is a rigid Lie algebra with trivial center whose radical is non nilpotent. This implies
that all derivations are inner. Let $f$ be a non trivial derivation of   $\mathcal{A}_q$. 
Since $\mathcal{A}_q$ is commutative, it is necessarily external.
Then $Id \otimes f$ is a derivation of $\frak{g}$ and satisfies 
$(Id \otimes f)(X \otimes 1)=X \otimes f(1)=0$ because $f(1 \cdot 1)=2f(1)=f(1)=0. $ Suppose that 
$Id \otimes f \in Int(\frak{g})$, that is $Id \otimes f=ad(\sum \alpha_{ij} X_i \otimes x_j).$ Thus
$(Id \otimes f)(X \otimes 1)=\sum \alpha_{ij}[X_i,X]\otimes x_j=0$ which implies $\sum \alpha_{ij}[X_i,X]=0$ for
any $j$ and $X.$ So $\sum \alpha_{ij}X_i \in Z(\frak{g}_p)$ for any $j.$ Since the center is
trivial, then $\sum \alpha_{ij} X_j=0$ for any $j$ and $Id \otimes f \notin Int(\frak{g}).$ There is a contradicion.
Therefore $\mathcal{A}_q$ is such that any external derivation is trivial. We deduce that $\mathcal{A}_q=M_1^q.$

\section{Cohomology and deformations}

\noindent a) The Chevalley cohomology of current Lie algebras was computed in \cite{Z}
for the degrees $1$ and $2$. It is shown that the algebra of derivation satisfies
$$\mathcal{D}er(\frak{g}\simeq \frak{g}_p \otimes \mathcal{A}_q )=\mathcal{D}er(\frak{g}_p)\otimes \mathcal{A}_q
\oplus Hom(\frak{g}_p/[\frak{g}_p,\frak{g}_p],Z(\frak{g}_p))]\otimes 
\frac{End(\mathcal{A}_q)}{\mathcal{A}_q+\mathcal{D}er \mathcal{A}_q}.$$
More precisely, let $f=f_1\otimes f_2$ be a derivation of $\frak{g}$. Denote $\mu _1$ the Lie product of 
$\g$ and $\mu _2$ the product of $\A _q$. Then $f$ is a derivation of $\g$ if and only if
$$\mu _1(f_1(X),Y)\otimes \mu _2(f_2(a),b)+\mu _1(X,f_1(Y))\otimes \mu _2(f_2(b),a)-f_1(\mu _1(X,Y))\otimes f_2(\mu _2(a,b))=0$$
for any $X,Y \in \g_1$ and $a,b \in \A _q.$ If $\mathcal{A}_q$ is unitary, by considering $a=b=1$, we obtain
$$[\mu _1(f_1(X),Y)+\mu _1(X,f_1(Y))-f_1(\mu _1(X,Y))]\otimes f_2(1)=0$$
and $f_1$ is a derivation of $\g_1$ as soon as $f_2(1) \neq 0$. Let us take $a=b$. The above identity
reduce to:
$$f_1(\mu _1(X,Y))\otimes (\mu _2(f_2(a),a)-f_2(a^2))=0.$$
Thus, either $f_1$ satisfies $f_1(\mu _1(\g_p,\g_p))=0$, or $f_2$ satisfies $\mu _2(f_2(a),a)=f_2(a^2).$ This last
identity becames  
$$\mu _2(f(a),b)+\mu _2(a,f(b))=2f_2(a,b)$$
by linearisation. Concerning the class of rigid Lie algebras that we consider, that is $Z(\frak{g}_p)=\{ 0 \}$, (note that
the conjecture that any rigid Lie algebra is of trivial center is still open), we deduce
$$\mathcal{D}er (\frak{g})\simeq \mathcal{D}er (\frak{g}_p)\otimes \mathcal{A}_q.$$

\noindent In the general case, the first space of cohomology is given in \cite{Z}:

$$H^1(\frak{g},\frak{g})\simeq H^1(\frak{g}_p,\frak{g}_p)\otimes \mathcal{A}_q \oplus Hom (\frak{g}_p,\frak{g}_p)\otimes
\mathcal{D}er (\mathcal{A}_q)\oplus Hom(\frak{g}_p/[\frak{g}_p,\frak{g}_p],Z(\frak{g}_p))]\otimes 
\frac{Hom(\mathcal{A}_q,\mathcal{A}_q)}{\mathcal{A}_q+\mathcal{D}er \mathcal{A}_q}.$$
In this case, this reduce to
$$H^1(\frak{g},\frak{g})\simeq H^1(\frak{g}_p,\frak{g}_p)\otimes \mathcal{A}_q \oplus Hom (\frak{g}_p,\frak{g}_p)\otimes
\mathcal{D}er (\mathcal{A}_q).$$
If $\frak{g}_p$ is rigid with non nilpotent nilradical, any derivation is inner. This implies
$H^1(\frak{g}_p,\frak{g}_p)=0$ and
$H^1(\frak{g},\frak{g})=Hom(\frak{g}_p, \frak{g}_p) \otimes \mathcal{D}er(\mathcal{A}_q).$ So 
$H^1(\frak{g},\frak{g})=0 \Leftrightarrow \mathcal{D}er(\frak{A}_q)=0.$
We find again the result.

\begin{proposition}
Let $\frak{g}_p$ be a rigid Lie algebra with a non nilpotent radical and a center reduced to zero. Then 
$\frak{g}=\frak{g}_p \otimes \mathcal{A}_q$ is rigid if and only if $\mathcal{D}er(\mathcal{A}_q)=0.$
\end{proposition}

\medskip

\noindent b) A Chevalley $2$-cochain $\varphi $ of $\frak{g}=\frak{g}_p \otimes \mathcal{A}_q$ decomposes as
$$\varphi =\psi _1\otimes  \varphi_2+\varphi _3\otimes \psi _4$$
with $\psi _1 \in \mathcal{C}^2 (\frak{g}_p,\frak{g}_p) $ , $\varphi _2 \in \mathcal{S}^2(\frak{g}_p,\frak{g}_p) $
and $\varphi_3 \in \mathcal{S}^2 (\frak{g}_p,\frak{g}_p)$,  $\psi _4 \in \mathcal{C}^2 (\mathcal{A}_p,\mathcal{A}_p) ,$
where $\mathcal{C}^2(\frak{g}_p,\frak{g}_p) $ denotes the space of Chevalley $2$-cochains of $\frak{g}_p$, 
$\mathcal{S}^2(\frak{g}_p,\frak{g}_p) $ the space of symmetric bilinear applications with values in 
$\frak{g}_p,$ $\mathcal{C}^2 (\mathcal{A}_p,\mathcal{A}_p)$ the space of
$2$-cochaines of the Harrison cohomology of $\mathcal{A}_q.$ 
We deduce using this decomposition that 
$H^2(\frak{g},\frak{g})=(H^2)' \oplus (H^2)''.$ The first space is compute in
(\cite{Z}, proposition 3.1). We find
$$(H^2)'=
H^2(\frak{g}_p,\frak{g}_p) \otimes \mathcal{A}_q \oplus \mathcal{B}(\frak{g}_p,\frak{g}_p) \otimes
\frac{H^2_H(\mathcal{A}_q,\mathcal{A}_q)}{\mathcal{P}_+(\mathcal{A}_q,\mathcal{A}_q)}
\oplus \chi (\frak{g}_p,\frak{g}_p) \otimes
\frac{\mathcal{A}(\mathcal{A}_q,\mathcal{A}_q)}{\mathcal{P}_+(\mathcal{A}_q,\mathcal{A}_q)}$$
(see  \cite{Z} for notations). But the second space was just computed when 
$\frak{g}_p$ is abelian. 

\noindent Let $\mu _1\otimes \mu _2 +\epsilon (\psi _1\otimes  \varphi_2+\varphi _3\otimes \psi _4)$ 
be an infinitesimal deformation of
$\mu _1\otimes \mu _2$.  The linear part of the Jacoby identity gives the expression of a $2$-cocycle of Chevalley 
cohomology of 
$\mu _1\otimes \mu _2$. We find:
$$
\begin{array}{ll}
\delta _{\mu _1\otimes \mu _2}(\psi _1\otimes  \varphi_2+\varphi _3\otimes \psi _4)= 
&\Sigma \mu _1(\psi _1(X_1,X_2),X_3)\otimes \mu _2(\varphi _2(a_1,a_2),a_3)\\
&+\Sigma \mu _1(\varphi _3(X_1,X_2),X_3)\otimes \mu _2(\psi _4(a_1,a_2),a_3)\\
& \Sigma \psi _1(\mu _1(X_1,X_2),X_3)\otimes \varphi  _2(\mu _2(a_1,a_2),a_3)\\
&+\Sigma \varphi_3(\mu _1(X_1,X_2),X_3)\otimes \psi _4(\mu _2(a_1,a_2),a_3)\\
& =0
\end{array}
$$
for any $X_1,X_2,X_3 \in \g_p$ and $a_1,a_2,a_3 \in \A_q$ and  the sum is taken on the cyclic permutations of 
$(1,2,3)$. 
We deduce
\begin{proposition}
If $\A_q$ is unitary then $\psi _1 \in Z^2(\g_p,g_p)$ as soon as $\varphi _2(1,1) \neq 0$.
\end{proposition}
If $X_1=X_2=X_3$, the above identity reduce to:
$$\mu _1(\varphi _3(X,X),X)\otimes \Sigma \mu _2(\psi _4(a_1,a_2),a_3)=0.$$
\begin{proposition}
If there exits $X \in \g_p$ such that $\mu _1(\varphi _3(X,X),X) \neq 0$ then $$\mu _2\bullet \psi _4=0$$
with
$$\mu _2\bullet \psi _4(a_1,a_2,a_3)=\Sigma \mu _2(\psi _4(a_1,a_2),a_3).$$
\end{proposition}
Note that $\psi $ is a $2$-cocyle for the Harrison cohomology of $\mu _2$ if 
$\mu _2\bullet \psi _4=\psi _4\bullet \mu _2.$

\noindent Suppose that $\g$ is rigid solvable with trivial center. Then $\A_q$ is unitay and 
$\psi _1 \in Z^2(\g_p,\g_p)$ as soon as
$\varphi _2(1,1) \neq 0.$

\section{Application : associative commutative real rigid algebras}

\subsection{Real rigid Lie algebras}
The study of the rigid real Lie algebras was lately initiated in 
\cite{A.C.G.G}. Let us point out the principal results. An external torus of derivations
of  $\frak{n}$ is an abelian subalgebra
$\frak{t}$ of
${\cal{D}}er(\frak{n})$, the Lie algebra of derivations of  $\frak{n}$, such as the elements are semi-simple.
This means that complex derivations $f \otimes Id \in \frak{t} \otimes \C$ are simultaneously diagonalizable.
If
$\frak{t}$ is a maximal (for the inclusion) external torus of $\frak{n}$ then
$\frak{t} \otimes \C$ is a maximal Malcev torus of $\frak{n} \otimes \C$. 
As all the maximal torus of
 $\frak{n} \otimes \C$ 
are conjugated with respect to 
 $Aut(\frak{n} \otimes \C)$, their dimensions are equal. It is the same for
the maximal torus 
 $\frak{t}$ of 
$\frak{n}$. This dimension is called the rank of $\frak{n}$. 
But contrary to the complex case, all the torus are not conjugated with respect to the group of automorphisms.
\begin{definition}
Let $\frak{n}$ be a finite dimensional real nilpotent Lie algebra. We call toroidal index of $\frak{n}$ 
the number of conjugaison classes of maximal external torus with respect to the group of aurtomorphisms
$Aut_{\R}(\frak{n})$ of $\frak{n}$.
\end{definition}
{\bf Example.} The toroidal index of the real abelian Lie algebra $\frak{a}_n$ of 
dimension $n$ is equal to $[n/2]+1$ where $[p]$ is the integer part of the rational number
$p$. In fact, let $\{X_1,...,X_n\}$ be a basis of  $\frak{a}_n$. Let us denote by
$f_i$ the derivation defined by $f_i(X_j)=\delta _i ^j X_j$ and by $f_{1,2p}$ the derivation given by
$$
\left\{
\begin{array}{l}
f_{1,2p}(X_{2p-1})=X_{2p},  \\
f_{1,2p}(X_{2p})=X_{2p-1}.\\
\end{array}
\right.
$$
Up a conjugation the maximal exterior torus are the subalgebras of $gl(n,\R)$ generated by 
$$
\begin{array}{ll}
{\frak{t}}_1=&\R\{f_1,...,f_n\} \\
{\frak{t}}_{2}=& \R\{f_{1,2},f_1+f_2,f_3...,f_n\} \\
{\frak{t}}_{3}=& \R\{f_{1,2},f_1+f_2,f_{1,4},f_3+f_4,f_5,...,f_n\} \\
... \\
{\frak{t}}_{n}=& \R\{f_{1,2},f_1+f_2,f_{1,4},f_3+f_4,...,f_{1,n},f_{n-1}+f_n\} \\
\end{array}
$$
if  $n$ is even, if not the last relation is replaced by 
$${\frak{t}}_{n}= \R\{f_{1,2},f_1+f_2,f_{1,4},f_3+f_4,...,f_{1,n-1},f_{n-2}+f_{n-1},f_n\}.$$

\subsection{Real rigid associative commutative algebras}
Let $\frak{r}_2$ be the real nonabelian $2$-dimensional Lie algebra. There exists a basis 
 $\{X_1,X_2\}$ with regard to which the bracket is given by $[X_1,X_2]=X_2.$ 
Let $\mathcal{A}_n$  be a $n$-dimensional real rigid commutative associative algebra.
Its complexified is isomorphic to
$M(1)_{\C}^n$. Thus the real current Lie algebra
$\frak{g}=\frak{r}_2 \otimes \mathcal{A}_n$ is rigid. We deduce that its complexified is rigid
and isomrphic to
$\frak{r}_2^n$. 
These remarks allow to write the following decomposition:
$$\frak{g}=\frak{r}_2 \otimes \mathcal{A}_n=\frak{t}_n \oplus \frak{a}_n$$
where $\frak{a}_n$ is the $n$-dimensional abelian Lie algebra.
We can deduct from this the structure of 
$\mathcal{A}_n$. 
In fact, if 
 $\{Y_1,...,Y_n\}$ is a basis of $\frak{t}_n$
corresponding to the derivations $f_{1,2},f_1+f_2,...,f_{1,2s},f_{2s-1}+f_{2s},f_{2s+1},...,f_n\}$
described in the previous section, the the Lie bracket
of  $\frak{g}$ satisfies
$$
\left\{
\begin{array}{l}
\lbrack Y_1,X_1]=-X_2, \ [Y_1,X_2]=X_1 \\
\lbrack Y_2,X_1]=X_1, \ [Y_2,X_2]=X_2 \\
... \\
\lbrack Y_{2s-1},X_{2s-1}]=-X_{2s}, \ [Y_{2s-1},X_{2s}]=X_{2s-1} \\
\lbrack Y_{2s},X_{2s-1}]=X_{2s-1}, \ [Y_{2s},X_{2s}]=X_{2s} \\
\lbrack Y_i,X_i]=X_i, \ i=2s+1,...,n.
\end{array}
\right.
$$
Let $\{e_1,...,e_n\}$ be a basis of $\mathcal{A}_n$ such that the isomorphism between
$\frak{r}_2 \otimes \mathcal{A}_n$ and $\frak{t}_n \oplus \frak{a}_n$ is given by
$U_1 \otimes e_i=Y_i$ and $X_{2i}=U_2 \otimes e_{2i-1}$, $X_{2i-1}=U_2 \otimes e_{2i}$
for $i=1,...,s$ and $X_j=U_2 \otimes e_j$ for $j=2s+1,..,n.$ 
The rigid associative algebra $\mathcal{A}_n$ is thus defined by 
$$
\left\{
\begin{array}{l}
e_{2i-1}^2=e_{2i-1}, \ \ i=1,...,s\\
e_{2i-1}e_{2i}=e_{2i }e_{2i-1}=e_{2i}, \ \ i=1,...,s\\
e_{2i}^2=-e_{2i-1}, \ \ i=1,...,s\\
e_j^2=e_j, \ \ j=2s+1,...,n.
\end{array}
\right.
$$
\begin{proposition}
Let $\mathcal{A}_n$ be a $n$-dimensional real rigid associative algebra. There exists an integer $s$, $1\leq s\leq n$
and a basis $\{e_1,...,e_n\}$ of $\mathcal{A}_n$ such that the multiplication of $\mathcal{A}_n$ is given
by 
$$
\left\{
\begin{array}{l}
e_{2i-1}^2=e_{2i-1}, \ \ i=1,...,s\\
e_{2i-1}e_{2i}=e_{2i }e_{2i-1}=e_{2i}, \ \ i=1,...,s\\
e_{2i}^2=-e_{2i-1}, \ \ i=1,...,s\\
e_j^2=e_j, \ \ j=2s+1,...,n.
\end{array}
\right.
$$
\end{proposition}

\end{document}